\documentclass[12pt]{article}

\usepackage{amsthm,amsmath,amsfonts,epsfig,amssymb,latexsym, bbm}
\usepackage{mathrsfs} 
\usepackage[T1]{fontenc}
\usepackage{graphicx}
\usepackage{enumerate}
\usepackage{xy,xypic}
\usepackage{graphics}
\usepackage{footnote}
\usepackage{scalefnt}
\usepackage{tikz}
\usepackage{color}
\usetikzlibrary{shapes,arrows,fit,calc,positioning}
\usetikzlibrary{matrix}
\xyoption{all}
\usepackage{imakeidx}
\makeindex

\title{Methods in complete intersections in corank one} 
\author{
 Satya Mandal \\ University of Kansas,  Lawrence, Kansas 66045, mandal@ku.edu\\
 \TCP{\bf ArrowticKTheory.com}
 }

\begin{document}

\pagenumbering{roman}
\setcounter{page}{0}

\renewcommand{\baselinestretch}{1.255}
\setlength{\parskip}{1ex plus0.5ex}
\date{15 August 2026} 
\newcommand{\iso}{\stackrel{\sim}{\longrightarrow}}
\newcommand{\sur}{\twoheadrightarrow}

\newcommand{\eop}{\hfill{\TCP{\rule{2mm}{2mm}}}}
\newcommand{\pf}{\noindent{\bf Proof.~}}
\newcommand{\outlinePf}{\noindent{\bf Outline of the Proof.~}}

\newcommand{\PD}{\dim_{{\SA}}}
\newcommand{\PDV}{\dim_{{\SV(X)}}}

\newcommand{\ra}{\rightarrow}
\newcommand{\lra}{\longrightarrow}
\newcommand{\hra}{\hookrightarrow} 
\newcommand{\Lra}{\Longrightarrow}
\newcommand{\Lla}{\Longleftarrow}
\newcommand{\Llra}{\Longleftrightarrow}
\newcommand{\pic}{The proof is complete.~}
\newcommand{\Sp}{\mathrm{Sp}}
\newcommand{\BiSp}{\mathrm{BiSp}}
\newcommand{\colim}{\mathrm{colim}}
\newcommand{\codim}{\mathrm{co}\dim}
\newcommand{\coh}{\mathrm{Coh}}
\newcommand{\qcoh}{\mathrm{QCoh}}

\newcommand{\bul}{\bullet}

\newcommand{\bE}{\begin{enumerate}}
\newcommand{\eE}{\end{enumerate}}

\newtheorem{theorem}{Theorem}[section]
\newtheorem{proposition}[theorem]{Proposition}
\newtheorem{lemma}[theorem]{Lemma}
\newtheorem{definition}[theorem]{Definition}
\newtheorem{corollary}[theorem]{Corollary}
\newtheorem{construction}[theorem]{Construction}
\newtheorem{notation}[theorem]{Notation}
\newtheorem{notations}[theorem]{Notations}
\newtheorem{remark}[theorem]{Remark}
\newtheorem{question}[theorem]{Question}
\newtheorem{example}[theorem]{Example} 
\newtheorem{examples}[theorem]{Examples} 
\newtheorem{exercise}[theorem]{Exercise} 
\newtheorem{setup}[theorem]{Setup} 
\newtheorem{openProblem}[theorem]{Open Problem} 
\newtheorem{Hyp}[theorem]{Hypothesis}

\newtheorem{problem}[theorem]{Problem} 
\newtheorem{conjecture}[theorem]{Conjecture} 

\newcommand{\bD}{\begin{definition}}
\newcommand{\eD}{\end{definition}}
\newcommand{\bP}{\begin{proposition}}
\newcommand{\eP}{\end{proposition}}
\newcommand{\bL}{\begin{lemma}}
\newcommand{\eL}{\end{lemma}}
\newcommand{\bT}{\begin{theorem}}
\newcommand{\eT}{\end{theorem}}
\newcommand{\bC}{\begin{corollary}}
\newcommand{\eC}{\end{corollary}} 
\newcommand{\bA}{\begin{application}} 
\newcommand{\eA}{\end{application}} 

\newcommand{\bN}{\begin{notation}} 
\newcommand{\eN}{\end{notation}} 

\newcommand{\bR}{\begin{remark}} 
\newcommand{\eR}{\end{remark}} 

\newcommand{\bExr}{\begin{exercise}} 
\newcommand{\eExr}{\end{exercise}} 

\newcommand{\bExm}{\begin{example}} 
\newcommand{\eExm}{\end{example}}

\newcommand{\TCP}{\textcolor{purple}}
\newcommand{\TCM}{\textcolor{magenta}}
\newcommand{\TCR}{\textcolor{red}}
\newcommand{\TCB}{\textcolor{blue}}
\newcommand{\TCG}{\textcolor{green}}
\def\spec#1{\mathrm{Spec}\left(#1\right)}
\def\proj#1{\mathrm{Proj}\left(#1\right)}
\def\supp#1{\mathrm{Supp}\left(#1\right)}
\def\Sym#1{{\CS}\mathrm{ym}\left(#1\right)} 
\def\norm#1{\parallel #1\parallel}
\def\LRf#1{\left( #1\right)}
\def\LRs#1{\left\{ #1\right\}}
\def\LRt#1{\left[ #1\right]}
\def\LBrace#1{\left\{\begin{array} #1\end{array}\right.}
\def\matrix#1{\left(\begin{array} #1\end{array}\right)}

\def\Hom#1{\underline{\mathrm{Hom}}\left(#1\right)}
\def\Obj#1{\underline{\mathrm{Obj}}\left(#1\right)}

\def\rank#1{\mathrm{rank}\left(#1\right)}

\def\Cat{\underline{\bf Cat}}
\def\eCat{\underline{\bf eCat}}
\def\dgCat{\underline{\bf dgCat}}

\def\a{\mathfrak {a}}
\def\b{\mathfrak {b}}
\def\c{\mathfrak {c}}

\def\d{\mathfrak {d}}
\def\e{\mathfrak {e}}
\def\f{\mathfrak {f}}
\def\g{\mathfrak {g}}
\def\i{\mathfrak {i}}
\def\j{\mathfrak {j}}
\def\k{\mathfrak {k}}
\def\l{\mathfrak {l}}
\def\m{\mathfrak {m}}
\def\n{\mathfrak {n}}
\def\p{\mathfrak {p}}
\def\q{\mathfrak {q}}
\def\r{\mathfrak {r}}
\def\s{\mathfrak {s}}
\def\t{\mathfrak {t}}
\def\u{\mathfrak {u}}
\def\v{\mathfrak {v}}

\def\w{\mathfrak {w}}

\def\x{\mathfrak {x}}
\def\y{\mathfrak {y}}
\def\z{\mathfrak {z}}

\def\A{\mathfrak {A}}
\def\B{\mathfrak {B}}
\def\C{\mathfrak {C}}
\def\D{\mathfrak {D}}
\def\E{\mathfrak {E}}

\def\F{\mathfrak {F}}

\def\G{\mathfrak {G}}
\def\H{\mathfrak {H}}
\def\I{\mathfrak {I}}
\def\J{\mathfrak {J}}
\def\K{\mathfrak {K}}
\def\L{\mathfrak {L}}
\def\M{\mathfrak {M}}
\def\N{\mathfrak {N}}
\def\O{\mathfrak {O}}
\def\P{\mathfrak {P}}

\def\Q{\mathfrak {Q}}

\def\R{\mathfrak {R}}

\def\Sf{\mathfrak {S}} 
\def\T{\mathfrak {T}}
\def\U{\mathfrak {U}}
\def\V{\mathfrak {V}}
\def\W{\mathfrak {W}}

\def\CA{\mathcal {A}}
\def\CB{\mathcal {B}}
\def\CP{\mathcal {P}}
\def\CC{\mathcal {C}}
\def\CD{\mathcal {D}}
\def\CE{\mathcal {E}}
\def\CF{\mathcal {F}}
\def\CG{\mathcal {G}}
\def\CH{\mathcal{H}}
\def\CI{\mathcal{I}}
\def\CJ{\mathcal{J}}
\def\CK{\mathcal{K}}
\def\CL{\mathcal{L}}
\def\CM{\mathcal{M}}
\def\CN{\mathcal{N}}
\def\CO{\mathcal{O}}
\def\CP{\mathcal{P}}
\def\CQ{\mathcal{Q}}

\def\CR{\mathcal{R}}

\def\CS{\mathcal{S}}
\def\CT{\mathcal{T}}
\def\CU{\mathcal{U}}
\def\CV{\mathcal{V}}
\def\CW{\mathcal{W}} 
\def\CX{\mathcal{X}}
\def\CY{\mathcal{Y}}
\def\CZ{\mathcal{Z}}

\newcommand{\smallcirc}[1]{\scalebox{#1}{$\circ$}}

\def\BA{\mathbb{A}}
\def\BB{\mathbb{B}}
\def\BC{\mathbb{C}}
\def\BD{\mathbb{D}}
\def\BE{\mathbb{E}}
\def\BF{\mathbb{F}}
\def\BG{\mathbb{G}}
\def\BH{\mathbb{H}}
\def\BI{\mathbb{I}}
\def\BJ{\mathbb{J}}
\def\BK{\mathbb{K}}
\def\BL{\mathbb{L}}
\def\BM{\mathbb{M}}
\def\BN{\mathbb{N}} 
\def\BO{\mathbb{O}}
\def\BP{\mathbb{P}} 
\def\BQ{\mathbb{Q}}
\def\BR{\mathbb{R}}
\def\BS{\mathbb{S}}
\def\BT{\mathbb{T}}
\def\BU{\mathbb{U}}
\def\BV{\mathbb{V}}
\def\BW{\mathbb{W}}
\def\BX{\mathbb{X}}
\def\BY{\mathbb{Y}}
\def\BZ{\mathbb {Z}}

\def\SA{\mathscr {A}} 
\def\SB{\mathscr {B}}
\def\SC{\mathscr {C}}
\def\SD{\mathscr {D}}
\def\SE{\mathscr {E}}
\def\SF{\mathscr {F}}
\def\SG{\mathscr {G}}
\def\SH{\mathscr {H}}
\def\SI{\mathscr {I}}
\def\SJ{\mathscr {J}}
\def\SK{\mathscr {K}} 
\def\SL{\mathscr {L}} 
\def\SM{\mathscr {M}}
\def\SN{\mathscr {N}}
\def\SO{\mathscr {O}}
\def\SP{\mathscr {P}}
\def\SQ{\mathscr {Q}}
\def\SR{\mathscr {R}}
\def\SS{\mathscr {S}}
\def\ST{\mathscr {T}}
\def\SU{\mathscr {U}}
\def\SV{\mathscr {V}}
\def\SW{\mathscr {W}}
\def\SX{\mathscr {X}}
\def\SY{\mathscr {Y}}
\def\SZ{\mathscr {Z}} 
 
\def\bfA{\bf A}
\def\bfB{\bf B}
\def\bfC{\bf C}
\def\bfD{\bf D}
\def\bfE{\bf E}
\def\bfF{\bf F}
\def\bfG{\bf G}
\def\bfH{\bf H}
\def\bfI{\bf I}
\def\bfJ{\bf J}
\def\bfK{\bf K}
\def\bfL{\bf L}
\def\bfM{\bf M}
\def\bfN{\bf N} 
\def\bfO{\bf O}
\def\bfP{\bf P} 
\def\bfQ{\bf Q}
\def\bfR{\bf R}
\def\bfS{\bf S}
\def\bfT{\bf T}
\def\bfU{\bf U}
\def\bfV{\bf V}
\def\bfW{\bf W}
\def\bfX{\bf X}
\def\bfY{\bf Y}
\def\bfZ{\bf Z}

\maketitle







\begin{abstract}
Let $A$ denote an affine algebra over an algebraically closed field ${\BF}$, with $\dim A=d\geq 3$. 
 Due to the recent exciting activities regarding projective $A$-modules $P$, with $\rank{P}=d-1$ (corank one), we try to explore the corank zero methods of N. Mohan Kumar and M. P. Murthy, in the complete intersections theory. 
  We hypothesize regarding cancellation properties projective modules of corank one, and derive some of analogues of the results in corank zero case.  We obtain some affirmative results, when $A$ is an affine algebra over  $\overline{{\BF}}_p$. 
\end{abstract} 
\pagenumbering{arabic}

\section{Introduction} 
Unless stated otherwise,  $A$ will always denote an affine algebra over an algebraically closed field $k$,
 with $\dim A=d\geq 3$. 
 %
  The category of finitely generated projective $A$-modules  will be denoted by ${\SP}\LRf{A}$.
 The subcategory of  $P\in {\SP}\LRf{A}$, with $\rank{P}=n$ will be denoted by ${\SP}_{n}\LRf{A}$. 
 For $P\in {\SP}\LRf{A}$, we also write 
 $$
 co\rank{P}=\dim A -\rank{P}.
 $$ 
  
For  $P\in {\SP}_d\LRf{A}$, N. Mohan Kumar and M. P. Murthy  investigated the question, whether the vanishing of the top Chern class $C^d\LRf{P}=0$ would detect that $P\cong Q\oplus A$. They considered a set of questions on complete intersections, that go hand in hand, with the $\rank{P}=d$ case
\cite{M99, MM98, M94, M88, MK85,  MK84, MkM82}.  
 In his seminal paper, Murthy   proved the final result \cite[Thm 3.7]{M94}. We provide an updated version of \cite[Thm 3.7]{M94},
as follows.
\bT\label{seminalThm}{\rm 
Let $A$ be a reduced affine algebra over an algebraically closed field ${k}$, with $\dim A=d\geq 2$. Let $P$ be a projective $A$-module with $\rank{P}=d$. Define a Chern class 
\begin{equation}\label{murthyChern}
c^d\LRf{P}=\sum_{r=0}^d\LRf{-1}^r\LRt{\Lambda^rP} \in F^dK_0\LRf{A} \subseteq K_0\LRf{A}.
\end{equation}
Then $P\cong Q \oplus A  \Llra c^d(P)=0$. 
}\eT
The Chern class theory makes best sense in the smooth case. While the definition of the top Chern class 
(\ref{murthyChern}) above, may appear ad-hoc, in smooth cases it coincides, via isomorphisms,  with the usual  top Chern class, which take value in the Chow group of zero cycles.

Analogous to the result (\ref{seminalThm}) in the $co\rank{P}=0$ case,   a celebrated recent result due to Asok, Bachmann and Hopkins \cite[Thm 3]{ABH26}, in the $co\rank{P}=1$ case, is the following.
\bT\label{ABH26Thm3}{\rm 
Let $X=\spec{A}$ be a  smooth affine scheme, over an algebraically closed field $k$ with $char\LRf{k}=0$, and $\dim X=d\geq 2$. Let $P$ be a projective $A$-module with $\rank{P}=d-1$. Then 
$$
P\cong Q\oplus A\quad \Llra \quad C^{d-1}\LRf{P}=0 \in CH^{d-1}\LRf{X}
$$ 
({\it Here and henceforth}, $CH^r\LRf{X}$ {\it will denote the Chow group of codimension $r$ cycles of} $X$, 
{\it and} $C^r(P)$ {\it will denote the $r^{th}$-Chern class of $P$. We further comment that we cannot replace 
$CH^{d-1}\LRf{X}$ by a subgroup of $K_0\LRf{X}$, in this statement}.)
}
\eT 
It would not be an overstatement that this result (\ref{ABH26Thm3}) has generated a considerable amount of excitements regarding the complete intersection questions in corank one case.

%
The work of Murthy \cite{M94}, comprehensively settled this project  on the entirety of the set of complete intersection questions, related to $P\in {\SP}_d\LRf{A}$, the corank zero case. ({\it By related questions we mean, for example,  whether or when  a locally complete intersection ideal $I\subseteq A$, with $height\LRf{I}=d$,  is  surjective image of a projective $A$-module $P\in {\SP}_d\LRf{A}$}.) Therefore, this article is about similar questions related to $P\in {\SP}_{d-1}\LRf{A}$ and local complete intersection ideals $J$ of height $d-1$.

  Technically,  the results on complete intersections depend, often enough, on the questions of whether  stably free modules, of certain rank, are free or not. By the cancellation theorem of Bass, if $\rank{P}\geq d+1$, then $P$ has the cancellation property. The following theorem of Suslin \cite[Thm 1]{S77b} played a key role in the proof of (\ref{seminalThm}) and related results. 
 \bT\label{SuslinStable}{\rm 
 Let $A$ be an affine algebra, over an algebraically closed field $k$, with $\dim A=d$. Let $P$ be a stably free $A$-module, with $\rank{P}=d$. Then $P$ is free. 
 }
 \eT
 In the same paper \cite{S77b}, Suslin asked what is the smallest integer $m$ such that, under the hypothesis 
 of (\ref{SuslinStable}), any stably free $A$-module $P$, with $\rank{P}\geq m$, is free? In response to this question, 
 N. Mohan Kumar \cite{MK85} constructed a  sequence of smooth affine algebras $A_p$, with $\dim A_p=p+2$,  where $p\geq 2$ is a prime number, 
 together with a stably free non-free  
 $A_p$-module $P_p$ with  $\rank{P_p}=p=\dim A_p-2$. 
 These examples 
 $P_p\in {\SP}_{p}\LRf{A_p}$,  may be considered as a 
 milestone in the study of Projective modules.  So, with respect to the question of Suslin, only case remains is 
 the case of $co\rank{P}=1$.
 We have the following from \cite{FRS12} on the $co\rank{P}=1$ case. 
 \bT\label{frsBara}{\rm 
 Let $A$ be a normal affine algebra over an algebraically closed field $k$, with $\dim A=d$ and $char(k)=p$. Assume $gcd\LRf{(d-1)!, p}=1$. If $d=3$, assume $A$ is smooth over $k$. Then every stably free $A$-module $P$, with $\rank{P}=d-1$, is free. 
 }
 \eT

In the spirit of  \cite{M94},  our contention in this article is that, the  vanishing theorems  (\ref{ABH26Thm3}) should be a consequence of certain stably free modules being free ({\it a cleaner version of theorem} \ref{frsBara}). 
 In the light of the approach  and results in \cite{M94}, following complete intersection statements make sense,
 which are associated to the case  $\rank{P}=d-1$. In the following statements, we assume $A$ is an affine algebra, over an algebraically closed field $k$, with $\dim A=d\geq 3$.

 \bE
 \item\label{canceFRS} {\bf (Stably Free modules:)}
 Is every stably free module $P\in {\SP}_{d-1}\LRf{A}$ free?
 \item\label{compInt}{\bf (Complete Intersections:)} Given a locally complete intersection ideal $J$ with 
 $height(J)=d-1$, when or whether there is a surjective map $\varphi: P\sur J$, with  $P\in {\SP}_{d-1}(A)$. 
 ({\it Note, such a surjective map $\varphi$ implies $\LRt{\frac{A}{J}}=\sum_{r=0}^{d-1}(-1)^r \LRt{\Lambda^rP}$ in $K_0(A)$. It also determines $C^{d-1}(P)=\LRt{cycle\LRf{\frac{A}{J}}}\in CH^{d-1}\LRf{A}$, in the smooth case.})
 \item \label{splingQn}{\bf (Splitting:)} Suppose $P\in {\SP}_{d-1}(A)$. Under what condition $P\cong Q\oplus A$ for some  $Q\in {\SP}_{d-2}(A)$? In particular, suppose there is a surjective map $\varphi: P\sur J$,  where 
 $J=\LRf{y_1, \ldots, y_{d-1}}$ is a complete intersection ideal, of height $d-1$. Does it imply $P\cong Q\oplus A$? 
 \item \label{chernQn}{\bf (Chern Class Vanishing:)} Assume $A$ is smooth over $k$. Let $P\in {\SP}_{d-1}(A)$ and the top Chern class $C^{d-1}(P)=0$. Does it follow $P\cong Q\oplus A$? 
 \eE 

 Contrary to the corank zero case (\ref{seminalThm}), in the corank one case, the vanishing Theorem \ref{ABH26Thm3} was not a consequence of the cancellation Theorem \ref{frsBara}. 
 Further, this pair of results (\ref{frsBara}, \ref{ABH26Thm3}) did not coalesce with 
 other natural questions on complete intersections. 
 Considering the existing theory on ${\SP}_d\LRf{A}$, it makes sense to expect that all four themes (\ref{canceFRS}, \ref{compInt}, \ref{splingQn}, \ref{chernQn}) mentioned above, and these  results (\ref{frsBara}, \ref{ABH26Thm3})
will 
grow into a comprehensive theory.

 In this article, we explored to implement the methods from the existing theory on ${\SP}_d\LRf{A}$ \cite{M94}, into the corank one case, 
  to obtain results on complete intersections and splitting (\ref{compInt}, \ref{splingQn}). 
 The standing normality hypothesis in (\ref{frsBara}), prevents the usual inductive argument to go through. 
 For a normal affine ring $A$, and a nonzero divisor $\t\in A$, the quotient ring $\frac{A}{A\t}$ is not normal, in general. %
  There is no evidence that the normality hypothesis is necessary in  (\ref{frsBara}).
 For our purpose this is not helpful and we are forced to write down a hypothesis, as follows.
 \begin{Hyp}[Hypothesis 1]\label{HypIntro}{\rm 
 Suppose $A$ is an affine algebra over an algebraically closed field and $\dim A=d\geq 3$. Given any 
 ideal ${\SI} \subseteq A$, denote the quotient algebra $\overline{A}=\frac{A}{\SI}$.  Then whenever 
 $\dim\LRf{\overline{A}}\leq d-1$,  any stably free $\overline{A}$-module $Q$, with $\rank{Q}=d-2$ is free. 
 }
 \end{Hyp}
Note that  the hypothesis is valid when $d=3$. Only other result known to us when the hypothesis (\ref{HypIntro}) is valid is the following theorem of Dhorajia and  Keshari \cite{DK12}.
\bT\label{Dkay12Thn}{\rm 
Let $p\geq 2$ be a prime number, and ${\BF}_p$ denote the prime field of order $p$. 
Let $A$ be an affine algebra over $\overline{{\BF}}_p$, with $\dim A=d \geq 4$. Assume $p\geq d$. Let $P\in {\SP}_{d-1}\LRf{A}$. Then $P$ has the cancellation property. 
}\eT
Under this hypothesis (\ref{HypIntro}), we prove a number of interesting consequences, analogous to the results of  Mohan Kumar and Murthy \cite{M94, MK84}. Note that Hypothesis \ref{HypIntro} imposes no restrictions on cancellation properties of projective $A$-modules. 

First, under (\ref{HypIntro}), we prove  an analogue of the key theorem in  
\cite[Thm 1.3]{M94}. Consequently, under (\ref{HypIntro}), we derive three  implications regarding complete intersections  (\ref{compInt})  and splitting (\ref{splingQn}).  In the following statements $I, J\subseteq A$ are ideals, with $I+J=A$, and  $Q\in {\SP}_{d-1}(A)$. We prove that Hypothesis \ref{HypIntro} implies the following:
\bE
\item\label{Aug15One} Suppose there is a  surjective map $Q\sur IJ$, and $J$ is a complete intersection ideal of height $d-1$. Then there is $F\in  {\SP}_{d-1}(A)$ with $\LRt{Q}=\LRt{F} \in K_0(A)$ and a surjective map $F \sur I$. 
Therefore, if modules in  ${\SP}_{d-1}(A)$ have the cancellation property, then there is a surjective map  $Q\sur I$.
\item\label{TwoAug15} Suppose  there is a surjective map $Q\sur J$, where  $J=\LRf{y_1, \ldots, y_{d-1}}$ is a complete intersection ideal of height $d-1$. Then $Q\oplus A\cong Q_0\oplus A^2$. Therefore, if modules in  ${\SP}_{d-1}(A)$ have cancellation property then $Q \cong Q_0\oplus A$. 
({\it Reminiscent to methods in} \cite{M94, MK84, M98s} {\it  this can be seen as an alternate version of 
Theorem \ref{ABH26Thm3}}). 
\item\label{3ree15Aug} Suppose $I$ and $J$ are locally complete intersection ideals of height $d-1$. If two of the three ideals $I$, $J$ and $IJ$ are images of stably free modules in ${\SP}_{d-1}(A)$ then so is the third one. 
\eE  

Since hypothesis \ref{HypIntro} holds, when $A$ is an affine algebra over $\overline{{\BF}}_p$, all the three statements 
(\ref{Aug15One}, \ref{TwoAug15}, \ref{3ree15Aug}) are valid for such a ring $A$, which we report in Section \ref{someRes}. Since these   concerns positive characteristic base field, there is no overlap between these  
affirmative results (Sec \ref{someRes}) and Theorem \ref{ABH26Thm3}. 

The key implications of hypothesis \ref{HypIntro} will be discussed in Section \ref{keyImpl}, and hypotheticals are discussed in Section \ref{SecHpo}. 
If and when Theorem \ref{frsBara} is improved and normality condition on $A$ is removed, 
(\ref{Aug15One}, \ref{TwoAug15}, \ref{3ree15Aug}) would be valid results, for all such rings. 
%


\section{The Key Implications}\label{keyImpl}
In this section, we derive the key implications of the Hypothesis  \ref{HypIntro}. 
We work under the following setup.
\begin{setup}\label{setUp}{\rm 
Unless specified otherwise,  $A$ will denote an  affine algebra over an algebraically closed field $k$, with $\dim A=d\geq 3$.
}
\end{setup}

\bP\label{asliDaDa}{\rm 
Let $A$ be as in (\ref{setUp}), satisfying the Hypotheses \ref{HypIntro}. Further assume that $A$ is reduced.  
Suppose $I$ and $J$ are two ideals with $I+J=A$. Assume $J$ is a locally complete intersection ideal 
with $height(J)=d-1$. Let  $P$, $Q$ be projective $A$-modules, and let
$$
\LBrace{{l}
\diagram P \ar@{->>}[r]^{\psi} & J\\ \enddiagram \\
\diagram Q \ar@{->>}[r]_{\varphi} & IJ\\ \enddiagram \\
}
\qquad\qquad  \LBrace{{l}{\rm be~surjective ~maps,}\\ 
{\rm where}~~P, Q\in {\SP}_{d-1}\LRf{A}. \\ 
{\rm Further,}~P~{\rm is~stably~free.}\\
}
$$  
Then there is a surjective map $\diagram P_0\ar@{->>}[r] & I \\ \enddiagram$, where $P_0\in {\SP}_{d-1}(A)$  with  $P_0\oplus P\cong Q\oplus A^{d-1}$. 
}
\eP
\pf Suppose that there is a surjective map $\diagram A^{d-1}\oplus Q \ar@{->>}[r]^{\f} & I \oplus P\\ \enddiagram$. 
Let $P_0={\f}^{-1}\LRf{I\times 0}$. Then $\diagram P_0\ar@{->>}[r] & I \\ \enddiagram$ is surjective. Further, 
$P_0\oplus P \cong Q\oplus A^{d-1}$. Therefore, we will construct such a map $\f$. Define the map
$$
\g: I \oplus P \lra A \quad by \quad \g(a, p)=\psi(p) -a.
$$ 
Then $\g$ is a surjective map. 
Consider the following  diagram of maps:  
$$
\diagram
0\ar[r] & \g^{-1}(IJ) \ar[r] \ar[d]_{\pi_0}^{\wr}&I\oplus P\ar[d]^{\pi} \ar@{->>}[r]^{\g} & A\\
0 \ar[r]&\psi^{-1}(IJ)\ar[r] &P \ar[r]_{\psi} & J \\
\enddiagram
$$ 
Here $\pi$ is the projection map.  Suppose $\g(a, p)=\psi(p)-a\in IJ$. Then $\psi(p)\in I \cap J=IJ$. Therefore, the 
restriction of $\pi$ to $\g^{-1}\LRf{IJ}$ defines the map $\pi_0$, as above. Suppose 
$\pi_0\LRf{a, p}=0$. Then $p=0$ and $\LRf{a, p}=\LRf{0,0}$. So, $\pi_0$ is injective. Let $p \in \psi^{-1}\LRf{IJ}$. 
Then $\psi(p) \in IJ$, and $\LRf{0, p}\in \g^{-1}\LRf{IJ}$. So, $\pi_0\LRf{0, p}=p$. So, $\pi_0$ is an isomorphism. 

Suppose that there is a surjective map $\chi: A^{d-2}\oplus Q \sur \psi^{-1}\LRf{IJ}$. Since $\g$ splits, such a map 
$\chi$ would lead to a surjective map $\diagram A^{d-1}\oplus Q \ar@{->>}[r]^{\f} & I \oplus P\\ \enddiagram$, as proposed above. So, we will define such a surjective map $\chi$. 
Now, consider the commutative diagram: 
\begin{equation}\label{HypoDiagmrnr}
\diagram 
0 \ar[r] & L \ar[r]\ar[d] & Q \ar@{-->}[dr]^{\eta_0} \ar[rr]^{\varphi}\ar[dd]^{\eta} && IJ \ar[r] \ar@{=}[d]& 0\\
&  K \ar[d]_1&&\psi^{-1}IJ\ar@{^(-->}[dl]\ar@{-->}[r]^{\psi_0} & IJ \ar[r]\ar@{^(->}[d] &0\\
0 \ar[r] & K \ar[r]  \ar[d]& P \ar[rr]_{\psi}  \ar[d]_{\p}&& J \ar[r] \ar[d]& 0\\
0 \ar[r] &\frac{K}{\eta L}\ar[d]\ar[r] &\frac{P}{\eta Q}\ar[d]\ar[rr]_{\beta} &&\frac{A}{I}\ar[d]\ar[r] & 0\\
 & 0 & 0& & 0 &\\
\enddiagram
~ \LBrace{{l}
{\rm The~
map}~\eta: Q\lra P~{\rm exists},  \\
{\rm 
because} ~Q~{\rm is~projective.}\\
{\rm Denote}, \\
L =\ker(\varphi), ~
K=\ker(\psi)\\
\psi_0~{\rm is ~the~restrictions~of}~\psi\\
\eta~{\rm factors~through}~\eta_0\\
{\rm The~other~maps~are~induced}.\\
}
\end{equation}
Since $J$ is locally complete intersection, 
it follows that the map $\overline{\eta}: \frac{Q}{JQ} \iso \frac{P}{JP}$, 
induced by $\eta$ is an isomorphism. So, $\eta_{1+s_0}$ is an isomorphism for some $s_0\in J$.  In fact, $P=\eta(Q)+JP$. 
By Nakayama's lemma there is $\t:= 1+s\in1+J$ such that $\t P \subseteq \eta Q$. 
Let ${\SI}=ann \LRf{\frac{K}{\eta L}}$. Since $\t P \subseteq \eta Q$, it follows
  that $\t\in {\SI}$ and ${\SI}+J=A$.  Further, $\frac{J}{{\SI}J}\cong \frac{A}{{\SI} A}$. 
By an application of the Snake lemma to the   $1^{st}$ and the (broken) $2^{nd}$ line of the diagram (\ref{HypoDiagmrnr}), we obtain the following
$$
  \frac{K}{\eta L} \cong 
\frac{\psi^{-1}IJ}{\eta Q}\qquad {\rm isomorphism}. 
$$ 
 Assume that there is a nonzero divisor $\t\in {\SI}$. So, with $\overline{A}=\frac{A}{{\SI}}$, $\dim \overline{A}\leq d-1$. 
Tensor the $3^{rd}$ line of (\ref{HypoDiagmrnr}) by $\overline{A}$. This leads to the exact sequence
$$
\diagram 
0 \ar[r] & \frac{K}{{\SI} K} \ar[r] & \frac{P}{{\SI} P} \ar[r] & \frac{A}{{\SI} A} \ar[r] & 0\\
\enddiagram
$$
It follows from Suslin's theorem (\ref{SuslinStable}) 
that $\frac{P}{{\SI}P} \cong \overline{A}^{d-1}$. By  hypothesis (\ref{HypIntro}), it follows that 
$\frac{K}{{\SI} K} \cong \overline{A}^{d-2}$.
By lifting this isomorphism to a surjective map  $\s: A^{d-2} \sur \frac{K}{{\SI} K}$, we construct   a  map $\theta$, as follows:
$$
\diagram 
&A^{d-2} \ar@{->>}[r]^{\s} \ar@/_/@{-->}[dl]_{\theta}\ar@{-->>}[d] & \frac{K}{{\SI} K} \ar@{->>}[d]\\
\psi^{-1}IJ \ar@{->>}[r]&\frac{\psi^{-1}IJ}{\eta Q} & \frac{K}{\eta L}\ar[l]_{\quad\sim} \\
\enddiagram
$$
So, the map
$$
\diagram 
A^{d-2} \oplus Q \ar@{->>}[rr]^{\chi:=\theta \oplus \eta_0} && \psi^{-1}IJ\\
\enddiagram 
$$ 
is surjective, as stipulated above. Therefore, the proof is complete, when we have a choice of $\t$, which is a nonzero divisor. 

In case, we do not have a choice of $\t\in {\SI}$, which is a nonzero divisor, we modify the map $\eta$, as follows. 
Let $\LRs{{\wp}_1, \ldots, {\wp}_r, {\wp}_{r+1}, \ldots , {\wp}_{r+r_1}}=\min{A}$ be the minimal primes of $A$. Assume 
$$
\LBrace{{l}
I \subseteq {\wp}_i ~\forall i=1, \ldots, r\\
I \not\subseteq {\wp}_i ~\forall i=r+1, \ldots, r+r_1\\
}
$$
Since $height(J) \geq 1$ we have $J \not\subseteq {\wp}_i~\forall i$. 
Subsequently, $\mu\LRf{-}$ will denote the number of generators of a module. Denote $N=K+\eta Q$. 
For $i=r+1, \ldots, r+r_1$, localize $1^{st}, 3^{rd}$ lines of (\ref{HypoDiagmrnr}), at $\wp={\wp}_i$:
$$
\diagram 
0 \ar[r] & L_{\wp} \ar[r]\ar[d] & Q_{\wp}  \ar[d]^{\eta} \ar[rr]^{\varphi}&& A_{\wp}  \ar[r] \ar@{=}[d]& 0\\
0 \ar[r] & K_{\wp}  \ar[r]  & P_{\wp}  \ar[rr]_{\psi}  && A_{\wp}  \ar[r] & 0\\
\enddiagram
 \Lra \LBrace{{l}
N_{\wp}=K_{\wp}+\eta \LRf{Q_{\wp}}=P_{\wp}.~So,\\
\mu\LRf{\LRf{\frac{P}{N}}_{\wp}}=0 \leq \mu\LRf{P_{\wp}}-(d-1)\\
}
$$
So, $N$ is $d-1$ fold basic at ${\wp}={\wp}_i$, for $i=r+1, \ldots, r+r_1$. 
For $i=1, \ldots , r$ and ${\wp}= {\wp}_i$, we have $I_{\wp}=0$. So, localization of $1^{s}, 3^{rd}$ lines of (\ref{HypoDiagmrnr}), leads to:
$$
\diagram 
0 \ar[r] & L_{\wp} \ar[r]^{\sim}\ar[d] & Q_{\wp}  \ar[d]^{\eta} \ar[rr]^{\varphi}&& I_{\wp} =0  \ar[d]& \\
0 \ar[r] & K_{\wp}  \ar[r]  & P_{\wp}  \ar[rr]_{\psi}  && A_{\wp}  \ar[r] & 0\\
\enddiagram
\Lra 
\LBrace{{l}
\mu\LRf{\LRf{\frac{P}{N}}_{\wp}}\leq \mu\LRf{\LRf{\frac{P}{K}}_{\wp}}= 1\\
=\mu\LRf{P_{\wp}} -(d-2) \\
}
$$
So, $N$ is $d-2$ fold basic at ${\wp}={\wp}_i$, for $i=1, \ldots, r$. 

By \cite[Lem 1.2]{M94}, there is a map $f\in Hom\LRf{Q, K}$ such that, with $\f:=\eta +f$, the image 
$\f\LRf{Q}$ is $d-1$ fold basic at ${\wp}_i$ for $i=r+1, \ldots, r+r_1$, and $d-2$ fold basic at 
 ${\wp}_i$ for $i=1, \ldots, r$. So,
 $$
 \LBrace{{ll}
 \mu\LRf{\LRf{\frac{P}{\f Q}}_{{\wp}_i}} \leq 0 & i=r+1, \ldots, r+r_1\\
  \mu\LRf{\LRf{\frac{P}{\f Q}}_{{\wp}_i}} \leq 1 & i=1, \ldots, r\\
 } 
 $$
 Therefore, if follows
 $$
  \LBrace{{ll}
\LRf{\frac{K}{\f L}}_{\wp_i}= \LRf{\frac{P}{\f Q}}_{\wp_i}=0   & i=r+1, \ldots, r+r_1\\
  \mu\LRf{\LRf{\frac{P}{\f Q}}_{{\wp}_i}} \leq 1 & i=1, \ldots, r\\
 } 
 $$
 So, for $i=r+1, \ldots , r+r_1$, the maps $\f_{\wp_i}$ are surjective, hence  isomorphisms. Therefore,
 $\LRf{\frac{K}{\f(L)}}_{\wp_i}=0$. 
 For $i=1, \ldots, r$ and ${\wp}:={\wp}_i$, consider the commutative diagram of exact sequences ({\it use  Snake lemma}):
 $$
\diagram 
0 \ar[r] & L_{\wp} \ar[r]^{\sim}\ar[d] & Q_{\wp}  \ar[d]^{\f} \ar[rr]^{\varphi}&& I_{\wp} =0  \ar[d]& \\
0 \ar[r] & K_{\wp}  \ar[r] \ar@{->>}[d] & P_{\wp} \ar@{->>}[d] \ar[rr]_{\psi}  && A_{\wp}  \ar[r] \ar[d]^1& 0\\
0\ar[r] &\LRf{\frac{K}{\f\LRf{L}}}_{\wp_i}\ar[r]&\LRf{\frac{P}{\f Q}}_{\wp_i}\ar[rr]_{\sim}&&A_{\wp_i}\ar[r]&0\\
\enddiagram 
$$
Considering, number of generators, in the local case, it follows $\LRf{\frac{K}{\f\LRf{L}}}_{\wp_i}=0$.

So, replacing $\eta$ by $\f$, if follows $\LRf{\frac{K}{\eta L}}_{\wp}=0~\forall {\wp}\in \min{A}$. Hence, there is nonzero divisor 
$\t\in ann\LRf{\frac{K}{\eta(L)}}$. So, the proof is complete by the previous case.
 $\eop$ 

\subsection{Hypotheticals} \label{SecHpo}
The proposition \ref{asliDaDa} is,  perhaps, the core of the algebraic arguments in \cite{MK84, M94}. 
Consequently,  (\ref{asliDaDa}) has a number of applications in 
 complete intersections, not withstanding our Hypothesis \ref{HypIntro}. The first one is motivated by \cite[Thm 2]{MK84}. 

\bC\label{BchiOne}{\rm 
Let $A$ be as in (\ref{setUp}), satisfying the Hypotheses \ref{HypIntro}. Further assume that $A$ is reduced. 
Let $J=\LRf{y_1, y_2, \ldots, y_{d-1}}$ be a complete intersection ideal, with $height(J)=d-1$. Suppose $I$ is an ideal, with $I+J=A$ and there is a surjective map $\diagram Q \ar@{->>}[r]^{\varphi} & IJ \enddiagram$, where $Q$ is a projective $A$-module with $\rank{Q}=d-1$. Then there is a surjective map $\diagram F \ar@{->>}[r]^{\f} & I \enddiagram$ where $F$ is a projective  $A$-module with $\rank{F}=d-1$, such that $\LRt{F}=\LRt{Q}\in K_0(A)$. 

}
\eC
\pf Follows from (\ref{asliDaDa}), with $P=A^{d-1}$. \pic $\eop$ 

\vspace{3mm}
The following is in the spirit of \cite[Cor 1.9]{M94}.
\bC\label{chuchiTwo}{\rm 
Let $A$ be as in (\ref{setUp}), satisfying the Hypotheses \ref{HypIntro}. Further assume that $A$ is reduced. 
Let $Q$ be a projective $A$-module with $\rank{Q}=d-1$.  Assume that there is a surjective map $\diagram Q \ar@{->>}[r]^{\f} & J \enddiagram$, where $J$ is a  complete intersection ideal of height $d-1$.
 Then $Q\oplus A\cong Q_0\oplus A^2$ for some $Q_0\in {\SP}_{d-2}\LRf{A}$.
}\eC
\pf Apply (\ref{asliDaDa}), with $I=A$, and $P=A^{d-1}$. Then there is a surjective map 
$\diagram P_1 \ar@{->>}[r]^{\f} & A \enddiagram$, for some projective $A$-module $P_1$, such that   
 $Q\oplus A^{d-1} \cong P_1\oplus A^{d-1}$. Clearly, $P_1=Q_0\oplus A$. Therefore, 
  $Q\oplus A^{d-1} \cong Q_0\oplus A^{d}$. By Suslin's theorem \cite{S77b, S77a}, we have $Q\oplus A\cong Q_0\oplus A^2$. 
  \pic $\eop$ 

\vspace{3mm}
The following is in the spirit of \cite[Cor 1]{MK84}.

\bC\label{evani3Ree}{\rm 
Let $A$ be as in (\ref{setUp}), satisfying the Hypotheses \ref{HypIntro}. Further assume that $A$ is reduced. 
Let $I$, $J$ be two locally complete intersection ideals with $height(I)=height(J)=d-1$, and $I+J=A$. If two of the three ideals $I$, $J$ and $IJ$ are images of  stably $A$-module of rank $d-1$, then so is the third one. 
}
\eC
\pf The case when $IJ$ is image of a stably free $A$-module, then the assertion follows directly from 
 (\ref{asliDaDa}).  Assume that there are surjective maps $\varphi_1: F_1 \sur I$ and 
 $\varphi_2: F_2 \sur J$ where $F_1$, $F_2$ are stably free $A$-modules of rank $d-1$. Then 
 $\varphi_1$ and $\varphi_2$ induce isomorphisms 
 $$
 \diagram
 \frac{F_1}{IF_1} \ar[r]^{\sim} & \frac{I}{I^2}\\
 \enddiagram 
 \quad  \diagram
 \frac{F_2}{JF_2} \ar[r]^{\sim} & \frac{J}{J^2}\\
 \enddiagram 
 $$
 Since $\dim\LRf{\frac{A}{I}}= \dim\LRf{\frac{A}{J}}=1$, both  $\frac{F_1}{JF_1}$, $ \frac{F_2}{JF_2}$
 are free. Combining, there is an isomorphism map $\overline{\varphi}: \frac{F}{IJF} \iso \frac{IJ}{I^2J^2}$,
 where $F$ is a free $A$-module, with $\rank{F}=d-1$. Now $\overline{\varphi}$ lifts to a surjective map
 $\varphi$, as follows
 $$
 \diagram 
 F\ar@{->>}[r]^{\varphi} \ar@{->>}[d]& IJK\ar@{^(->}[d]\\
  \frac{F}{IJF} \ar[r]_{\overline{\varphi}} &\frac{IJ}{I^2J^2}\\
 \enddiagram
 \quad {\rm where} \LBrace{{l} 
 K~{\rm is~a~local~complete ~intersection}\\
 {\rm ideal~with }~height(K)=d-1\\
\ni \quad IJ+K=A\\
 }
 $$
 Since $I$ is image of $F_1$, and $I(JK)$ is image of $F$, by  (\ref{asliDaDa}) there is surjective map 
 $F_3 \sur JK$, where $F_3$ is a stably free $A$-module with $\rank{F_3}=d-1$. By same argument,
 there is a surjective map $F_4 \sur K$ where $F_4$ is a stably free $A$-module with $\rank{F_4}=d-1$.
 Again by the same argument, we have a surjective map $F_5 \sur IJ$, where
 $F_5$ is a stably free $A$-module with $\rank{F_5}=d-1$.
 \pic $\eop$

\section{Affirmative results}\label{someRes}
For prime numbers $p\geq 2$, ${\BF}_p$ will denote the prime field of order $p$, and  $\overline{{\BF}}_p$ denotes the algebraic closure. 
Suppose $A$ is a reduced affine algebra over $\overline{{\BF}}_p$, with $\dim A=d\geq 3$. 
Then the Hypothesis \ref{HypIntro} is valid (Thm \ref{Dkay12Thn}), for such algebras, when $p\geq d$. Accordingly, the hypotheticals
in Section \ref{SecHpo}, are valid results in this case.
\bC\label{UboEin}{\rm 
Let $A$ be reduced affine algebra over $\overline{{\BF}}_p$, with $4 \leq \dim A=d\leq p$. 
Let $J=\LRf{y_1, y_2, \ldots, y_{d-1}}$ be a complete intersection ideal, with $height(J)=d-1$. Suppose $I$ is an ideal, with $I+J=A$ and there is a surjective map $\diagram Q \ar@{->>}[r]^{\varphi} & IJ \enddiagram$, where $Q$ is a projective $A$-module with $\rank{Q}=d-1$. Then there is a surjective map $\diagram Q \ar@{->>}[r]^{\f} & I \enddiagram$. %

}
\eC
\pf By theorem \ref{Dkay12Thn}, hypothesis \ref{HypIntro} is valid for $A$. By (\ref{BchiOne}), there is a surjective map $\diagram F \ar@{->>}[r] & I \enddiagram$ where $F\in {\SP}_{d-1}\LRf{A}$ and 
$\LRt{Q}=\LRt{F}\in K_0\LRf{A}$. Again, by (\ref{Dkay12Thn}), $Q\cong F$. \pic 
$\eop$

\vspace{3mm}
The following a counter part of the celebrated vanishing theorem of Asok-Bachman-Hopkins (\ref{ABH26Thm3}). 
%
The following is in the spirit of \cite[Cor 1.9]{M94}.
\bC\label{dno836}{\rm 
Let $A$ be reduced affine algebra over $\overline{{\BF}}_p$, with $4 \leq \dim A=d\leq p$. 
Let $Q\in {\SP}_{d-1}\LRf{A}$, such that there is  a surjective map $\diagram Q \ar@{->>}[r]^{\f} & J \enddiagram$, where $J$ is a  complete intersection ideal of height $d-1$.
 Then $Q\cong Q_0\oplus A$ for some $Q_0\in {\SP}_{d-2}\LRf{A}$.
}\eC
\pf By theorem \ref{Dkay12Thn}, hypothesis \ref{HypIntro} is valid for $A$.
It follows from (\ref{chuchiTwo}) that $Q\oplus A\cong Q_0\oplus A^2$, for some $Q_0\in {\SP}_{d-2}\LRf{A}$. 
By theorem \ref{Dkay12Thn}, $Q\cong Q_0\oplus A$. \pic $\eop$ 

\vspace{3mm}
The following is in the spirit of \cite[Cor 1]{MK84}.

\bC\label{Mohan333}{\rm 
Let $A$ be reduced affine algebra over $\overline{{\BF}}_p$, with $4 \leq \dim A=d\leq p$. 
Let $I$, $J$ be two locally complete intersection ideals with $height(I)=height(J)=d-1$, and $I+J=A$. If two of the three ideals $I$, $J$ and $IJ$ are complete intersections,  then so is the third one. 
}
\eC
\pf Follows form (\ref{evani3Ree}) and the cancellation theorem \ref{Dkay12Thn}. \pic $\eop$

\bR{\rm 
We close this discussions, with a few remarks.
\bE
\item 
The height $d$ version of (\ref{Mohan333}) was proved in \cite{MK84}, for any affine algebra over any algebraically closed field ${\BF}$. 
\item Suppose $A$ is an affine algebra  over any algebraically closed field ${\BF}$, with $\dim A=3$. Then the hypothesis \ref{HypIntro} is valid. Therefore, all the statements in Hypothetical section \ref{SecHpo} are valid results.
However, there is no cancellation theorem for $P\in {\SP}_2\LRf{A}$, unless $A$ is smooth. The case of smooth three folds was dealt with in \cite{MkM82}.
\item The results (\ref{UboEin}, \ref{dno836}, \ref{Mohan333}) in this section, apply only to the affine algebras $A$, over fields ${\BF}$ of positive characteristic $p\neq 0$. Consequently, there is no overlap with the Theorem \ref{ABH26Thm3}.
\eE

}\eR 






\begin{thebibliography}{200}


\bibitem[ABH26]{ABH26}  Aravind Asok, Tom Bachmann, Michael J. Hopkins
On $\BP^1$
-stabilization in unstable motivic homotopy theory, Annals of Mathematics (to appear).

%


%

%




%


\bibitem[DK12]{DK12} Dhorajia, Alpesh M.; Keshari, Manoj K
A note on cancellation of projective modules
J. Pure Appl. Algebra 216 (2012), no. 1, 126–129.


\bibitem[FRS12]{FRS12}
Fasel, J.; Rao, R. A.; Swan, R. G.
On stably free modules over affine algebras
%
Publ. Math. Inst. Hautes Études Sci. 116 (2012), 223–243.


%

%



 

 
 
 
 

















 










%



%

\bibitem[M98s]{M98s} Mandal, Satya
Complete intersection  K -theory and Chern classes
%
Math. Z. 227 (1998), no. 3, 423–454.

\bibitem[MM98]{MM98} Mandal, Satya; Pavaman Murthy, M.
%
Ideals as sections of projective modules
J. Ramanujan Math. Soc. 13 (1998), no. 1, 51–62.
%

 
 
 
 

 
 
 
 

 

 

 
 

 
 
 
 
 


\bibitem[MK85]{MK85} N Mohan Kumar, Stably Free Modules, 
American Journal of Mathematics , Dec., 1985, Vol. 107, No. 6 (Dec., 1985), pp.
1439-1444

\bibitem[MK84]{MK84} Mohan Kumar, N.
Some theorems on generation of ideals in affine algebras
%
Comment. Math. Helv. 59 (1984), no. 2, 243–252.

\bibitem[MkM82]{MkM82} Kumar, N. Mohan; Murthy, M. Pavaman
Algebraic cycles and vector bundles over affine three-folds
%
Ann. of Math. (2) 116 (1982), no. 3, 579–591.


\bibitem[M99]{M99} Murthy, M. Pavaman
%
A survey of obstruction theory for projective modules of top rank
Contemp. Math., 243
American Mathematical Society, Providence, RI, 1999, 153–174.

\bibitem[M94]{M94} Murthy, M. Pavaman
Zero cycles and projective modules
%
Ann. of Math. (2) 140 (1994), no. 2, 405–434.

\bibitem[M88]{M88} Murthy, M. Pavaman
 Zero-cycles, splitting of projective modules and number of generators of a module
%
Bull. Amer. Math. Soc. (N.S.) 19 (1988), no. 1, 315–317.













   
 
  
  

  
  
   

   %











 


  
   %
   

\bibitem[S77a]{S77a} Suslin, A. A.
A cancellation theorem for projective modules over algebras
%
Dokl. Akad. Nauk SSSR 236 (1977), no. 4, 808–811.

\bibitem[S77b]{S77b} Suslin, A. A.
%
Stably free modules
Mat. Sb. (N.S.) 102(144) (1977), no. 4, 537–550, 632.







%


%




%


 \end{thebibliography}
\end{document}